\documentclass[preprint,3p,11pt]{elsarticle}

\usepackage{bm}
\usepackage{amsmath}
 \usepackage{hyperref}
 \hypersetup{backref, colorlinks=true, linkcolor=blue}

\usepackage{stfloats}
\usepackage{subcaption}
\usepackage{placeins}

\usepackage{graphicx} 
\usepackage{booktabs} 
\usepackage{amssymb,amsmath}
\usepackage[ruled,vlined]{algorithm2e}
\usepackage{xcolor}

\begin{document}
\begin{center}
{\bf \Large Second-Moment  Method for Transport Problems  	with Anisotropic  \\ 
	\vspace{0.1cm}
Scattering}
\end{center}
\begin{center}
India J. Allan and Dmitriy Y. Anistratov 
\end{center}
\begin{center}
{ \it Department of Nuclear Engineering,
North Carolina State University,
Raleigh, NC 27695 \\
 ijallan@ncsu.edu, anistratov@ncsu.edu }
\end{center}

\vspace{-2.1cm}
\begin{frontmatter}
\begin{abstract}
This paper presents  a new nonlinear two-level  acceleration method for solving the particle transport equation   with anisotropic scattering. 
The method is formulated with the projection operator approach. 
The low-order equations  are defined for  the angular moments using   projection operators and closures of the second-moment  method. 
A nonlinear prolongation operator is applied  to the scattering term to derive the high-order transport equation.
The nonlinear system of high-order and low-order equations is equivalent to the original transport equation.
The   equations  are  approximated in space by the lumped linear-discontinuous   Galerkin method. 
Numerical results are presented to demonstrate the performance of the proposed numerical method.
\end{abstract}
\end{frontmatter}
\vspace{-1cm}

\section{Introduction}
Particle transport problems with highly anisotropic scattering arise 
in various applications such as radiation shielding, medical physics, 
atmospheric sciences, plasma physics, and astrophysics.
To solve this type of transport problems,   acceleration methods have been developed and theoretically analyzed
\cite{dya-vyag-kiam-1985,khattab=larsen-1991,morel-manteuffel-1991,vyag-aristova-m&c1997,pautz-morel-adams-m&c1999,aristova-vyag-2000,patel-warsa-prinja-2020}.
In this paper, we present a new transport acceleration method for the Boltzmann transport equation (BTE) with anisotropic scattering. 
We apply a projection operator approach to formulate a two-level nonlinear iteration method. 
The low-order equations for  the angular moments are derived  using   projection operators and closures of the second-moment (SM) method \cite{sm-1976}. 
We formulate the high-order transport equation by applying a nonlinear prolongation
operator \cite{dya-vyag-kiam-1985,vyag-aristova-m&c1997}.
The resulting nonlinear system of high-order and low-order equations is equivalent to the original BTE.
The system of equations is  discretized by the lumped linear-discontinuous (LLD) Galerkin method. 
We derive the discrete form of the nonlinear prolongation  operator for the high-order transport equation approximated with the LLD scheme.
Numerical results are presented to demonstrate the performance of the proposed computational method.

\section{Nonlinear Second-Moment Method}
The steady-state one-group   Boltzmann transport equation with anisotropic scattering in
1D slab geometry is given by
\begin{equation}    \label{t-eq}
	\mu 	\frac{\partial}{\partial x }  \psi(x, \mu)   + \sigma_{t}(x)  \psi(x, \mu)  =   
	\frac{\sigma_{s} }{4 \pi}
	\int_0^{2 \pi} d \gamma'	\int_{-1} ^1  d \mu'
	f_{s} (\mu_0)   \psi(x, \mu')   +  \frac{1}{4 \pi}q(x) \, ,
\end{equation}
\begin{equation*}
	x\in[0,X] \, ,\quad  \mu \in [-1,1]\, ,
\end{equation*}
\begin{align} \label{t-eq-bcs}
	\psi(0,\mu) = \psi_{in}^+(\mu)  \ \mbox{for} \ \mu>0\, ,\ \  \psi(X,\mu) = \psi_{in}^-(\mu) \ \mbox{for} \ \mu<0 \, .
\end{align}
Here $\psi$ is the angular flux; $x$ is the spatial position; $\mu$ is the directional cosine of particle motion;
$\gamma$ is the azimuthal angle;
$\sigma_t$ and $\sigma_s$  are the total and scattering  cross sections, respectively;
$f_s(\mu_0)$ is the  scattering probability density function (aka the phase function);
$\mu_0$ is the cosine of the scattering angle;  $q$ is the  isotropic source.

The low-order SM (LOSM) equations are derived by operating on the transport equation
with $2\pi \int_{-1}^{1} (\cdot) d \mu$ and $2 \pi \int_{-1}^{1} (\cdot) \mu d \mu$
as well as defining the closure for the highest angular moment in the projected equations given by \cite{sm-1976}
\begin{equation}
	P(x)= \frac{2\pi}{3} \int_{-1}^1 (1 -3 \mu^2 )  \psi(x,\mu) d \mu  \, .
\end{equation}
The LOSM for the scalar flux, $\phi=2\pi \int_{-1}^{1} \psi d \mu$
and current, $J=2\pi \int_{-1}^{1}  \mu\psi d \mu$ are defined  by
\begin{subequations} \label{losm-eqs}
	\begin{equation}   \label{losm-eq-0}
		\frac{d J}{dx}(x)  + \big(\sigma_t(x)- \sigma_{s,0}(x) \big) \phi (x)= q (x) \, ,
	\end{equation}
	\begin{equation} \label{losm-eq-1}
		\frac{1}{3}\frac{d \phi }{dx}(x)  + \big(\sigma_t (x)- \sigma_{s,1}(x)  \big) J(x)= 	\frac{d P}{dx}(x)    \, ,
	\end{equation}
	\begin{equation} \label{sm-bc-0}
		J(0)= -\frac{1}{2}\phi(0)   +2 J_{in}^+ + B_0 \, ,
			\end{equation}
	 \begin{equation} \label{sm-bc-X}
		J(X)= \frac{1}{2}\phi(X)   + 2 J_{in}^- -  B_X\, ,
	\end{equation}
\end{subequations}
where $ \sigma_{s,\ell} =  \sigma_{s} f_{s,\ell}$ and 
\begin{equation}
	f_{s,\ell}  = 2 \pi \int_{-1}^1  f_{s} (\mu_0) P_{\ell}(\mu_0) d \mu_0 \, 
\end{equation}
are the coefficients of expansion of the phase function in the Legendre polynomials  $P_{\ell}(\mu_0)$,
\begin{equation}
	J_{in}^{\pm} = \pm 2 \pi  \int_{0}^{\pm 1}  \mu \psi_{in}^{\pm} d \mu \,.
\end{equation}
The boundary closures are given by
\begin{subequations}\label{sm-boundary}
	\begin{equation}
		B_0=  2\pi \int_{-1}^{1} \bigg( \frac{1}{2} - |\mu| \bigg)  \psi(0,\mu) d\mu \, ,
	\end{equation}
	\begin{equation}
		B_X=   2 \pi \int_{-1}^{1} \bigg(\frac{1}{2} - |\mu|  \bigg)  \psi(X,\mu) d\mu \, .
	\end{equation}
\end{subequations}

To formulate a high-order transport equation,  
we  apply  a nonlinear prolongation operator 
which transforms  the scattering term into a form in terms of the scalar flux, current,
and  linear-fractional factors \cite{dya-vyag-kiam-1985,vyag-aristova-m&c1997}.
First, we cast the phase function using the first three terms of its Legendre polynomial expansion to obtain
\begin{equation} \label{fs-expansion}
	f_{s} (\mu_0)   =\frac{1}{4 \pi } \Big(    f_{s,0} + 3   f_{s,1} \mu_0 +  \frac{5}{2}   f_{s,2} (3 \mu_0^2 -1) \Big) +   f_s^* (\mu_0) \, ,
\end{equation}
where  $f_s^* (\mu_0)$ is the residual term of the expansion. We then introduce 
the Eddington factor 
\begin{equation} \label{edd-f}
	E(x) =   \frac{\int_{-1}^1  \mu^2\psi(x,\mu)d \mu }{\int_{-1}^1  \psi(x,\mu)d \mu} \, ,
\end{equation}
and the factor defined with the residual expansion term 
\begin{equation} \label{f-star-f}
	F_s^*(x,\mu) =	\frac{\int_{-1} ^1      f_s^* (\mu', \mu)    \psi(x, \mu') d \mu'}{\int_{-1}^1 \psi(x, \mu') d \mu'} \, ,
\end{equation}
where
\begin{equation}
	f_s^* (\mu', \mu)  =  \frac{1}{2\pi}\int_0^{2 \pi} d \gamma  \int_0^{2 \pi} d \gamma' 	   f_s^* (\mu_0) 
\end{equation}
and
\begin{equation}
	\mu_0=\cos (\gamma - \gamma')  \sqrt{(1- \mu^2)(1-(\mu')^2)} + \mu \mu' \, .
\end{equation}
We plug Eq. \eqref{fs-expansion} into the right-hand side  (RHS) of Eq.  \eqref{t-eq} and use factors  \eqref{edd-f} and  \eqref{f-star-f} to 
obtain the transport equation in the following form:
\begin{multline} \label{bte-closed}
	\mu 	\frac{\partial}{\partial x }  \psi(x, \mu)   +   \sigma_{t}(x)  \psi(x, \mu)  =     
	\frac{1 }{2}\sigma_{s}(x)
	\bigg(   f_{s,0} (x) 
	+ \frac{5 } {4}   f_{s,2}(x) (3\mu^2 -1) (3E(x) - 1) + 2 F_s^*(x,\mu)  \bigg)\phi(x)    \\
	+	\frac{3 }{2}  \sigma_{s}(x)f_{s,1}(x)  \mu J(x)     + \frac{1}{4 \pi}q(x)\, .
\end{multline}
The two-level SM method for the transport equation with anisotropic scattering is defined by the high-order transport equation \eqref{bte-closed} and LOSM equations
\eqref{losm-eqs}. 
This is a nonlinear method due to the nonlinearity of the RHS of Eq. \eqref{bte-closed}.
Hereafter we refer to this method as the nonlinear SM (NSM) method.
The system of equations in the NSM method is solved by the fixed-point iteration method. 
Algorithm \ref{NSM-iter} presents the  NSM iteration scheme, where $\phi_{in}^{\pm} = \pm \int_0^{\pm 1} \psi_{in}^{\pm}d \mu$.
\begin{algorithm}
	\DontPrintSemicolon
	$s=-1$\;
	$P^{(0)}\! = \! 0$, $B_0^{(0)}\! =\! \frac{1}{2} \phi_{in}^+ \!- \! J_{in}^+$,   $B_X ^{(0)}\! =\! J_{in}^- \! +  \! \frac{1}{2} \phi_{in}^-$,  $\rho^{(0)}\! = \! \rho^{(1)}\! =\! \frac{1}{2}$\;
	\While{$\| \phi^{(s)} -  \phi^{(s-1)} \|_{\infty} > \varepsilon \Big( \frac{1}{\rho^{(s)}}  - 1   \Big)$}{
		$s = s+ 1$\;
		\If{$s \ge 1$}{
			solve the high-order transport equation \eqref{bte-closed} for $\psi^{(s)}$ using 
			$\phi^{(s-1)}$, $J^{(s-1)}$, $E^{(s-1)}$, $F_s^{ \star (s-1)}$ \;
			compute $P^{(s)}$, $B_0^{(s)}$, $B_X^{(s)}$,  $E^{(s)}$, $F_s^{ \star (s)}$  using  $\psi^{(s)}$ \;
		}
		solve the LOSM equations \eqref{losm-eqs}  defined with $P^{(s)}$, $B_0^{(s)}$, and $B_X^{(s)}$ to compute $\phi^{(s)}$ and $J^{(s)}$\;
		\If{$s \ge 2$}{
			$\rho^{(s)} = \frac{ \|  \phi^{(s)} -  \phi^{(s-1)} \|_{\infty } } { \| \phi^{(s-1)} -  \phi^{(s-2)} \|_{\infty } }$
		}
	}
	\caption{The NSM method for anisotropic transport problems \label{NSM-iter}} 
\end{algorithm}

\section{Discretization of Equations}
We define the spatial grid $\{\omega_i \}_{i=1}^{n_x}$, where $\omega_i = [x_{i-1},x_i]$,
$x_0=0$ and $x_{n_x} = X$.
The cross sections are  piece-wise constant functions on the set of mesh cells.
To discretize the high-order transport equation  \eqref{bte-closed},
we apply the method of discrete ordinates and use the quadrature set $\{ \mu_m , w_m\}_{m=1}^{n_{\mu}}$.
The transport equation is discretized in space by the LLD Galerkin method. 
The linear discontinuous (LD) expansion of   the solution $\psi_m(x) = \psi (x, \mu_m)$ in space is given by
\begin{equation}\label{dfem-x-lagr} 
	\psi_{m,i}(x) = \sum_{\alpha=L,R} \psi_{m,i,\alpha} b_{i,\alpha}(x) \, , \quad
	x \in \omega_i \, , 
\end{equation}
where the basis functions are defined by
\begin{equation}\label{dfem-basis}
	b_{L,i} = \frac{1}{\Delta x_i} ( x_i - x) \, , \quad
	b_{R,i} = \frac{1}{\Delta x_i} ( x - x_{i-1}) \, , \quad \Delta x_i = x_i -x_{i-1} \, .
\end{equation}
The high-order transport equation \eqref{bte-closed} discretized with  the LLD method is given by
\begin{multline}  \label{bte-lld}
	\mu_m\Delta x_i \Big(  \mathbb{L}^b	\boldsymbol{\psi}_{m,i}^b     
	+  \mathbb{L} \boldsymbol{\psi}_{m,i}   \Big)   
	+ 	\Delta x_i     \sigma_{t,i}   \Big)    \mathbb{M}  \boldsymbol{\psi}_{m,i}  =   
	\frac{\Delta x_i}{2} \sigma_{s,i} \bigg[ \ f_{s,0,i}   \mathbb{M}   \boldsymbol{\phi}_{i}  
	+ 		  3  f_{s,1,i}  \mu_m \mathbb{M}   \boldsymbol{J}_{i}     \\
	+    \frac{5}{4}  f_{s,2,i}  (3\mu_m^2 -1)  \mathbb{M}\big(  3\boldsymbol{E}_{i} \odot \boldsymbol{\phi}_{i}   - \boldsymbol{\phi}_{i}   \big)  
	+  2    \mathbb{M}\boldsymbol{ { F}}_{s,m,i}^* \odot  \boldsymbol{\phi}_{i}     
	\bigg]
	+  \frac{1}{4 \pi} \Delta x_i  \mathbb{M}   \bm{q}_{i}     \, ,
\end{multline}
where $\odot$ is the Hadamard product,
\begin{equation} \label{L-matrices}
	\mathbb{L}^b =   \left[\begin{array}{rl} -1 & 0  \\ 0 & 1  \end{array} \right] \, ,
	\  \
	\mathbb{L} =   \frac{1}{2}\left[\begin{array}{rr}     1 & 1 \\ -1 & -1   \end{array} \right]  \, ,
	\  \
	\mathbb{M} =  \frac{1}{2} \left[\begin{array}{rl}   1 & 0  \\ 0 & 1   \end{array} \right] \, ,
\end{equation}
\begin{equation}
	\boldsymbol{\psi}_{m,i}  = [\psi_{m,i,L}, \psi_{m,i,R}  ]^{\top} \, ,
	\quad
	\boldsymbol{\psi}_{m,i}^b = [\psi_{m,i,L}^b, \psi_{m,i,R}^b  ]^{\top} \, ,
\end{equation}
\begin{equation}
	\psi_{m,i,L}^b = 
	\begin{cases}
		\psi_{in,m}^{+}   & \mu>0 \, , \, i=0 \, , \\
		\psi_{m,i-1,R}      & \mu>0 \, , \, i>0 \, , \\
		\psi_{m,i,L}     & \mu<0 \, ,
	\end{cases}
\end{equation}
\begin{equation}
	\psi_{m,i,R}^b = 
	\begin{cases}
		\psi_{in,m}^{-}  & \mu<0 \, , \, i=n_x \, , \\
		\psi_{m,i+1,L}   & \mu<0 \, , \, i<n_x \, , \\
		\psi_{m,i,R}        & \mu>0 \, ,
	\end{cases}
\end{equation}
\begin{equation}
	\boldsymbol{E}_{i}  = [E_{i,L}, E_{i,R}]^{\top} \, , \quad 
	\boldsymbol{ {F}}_{s,m,i}^*  = [F_{s,m,i,L}^* , F_{s,m,i,R}^*]^{\top}  \,   ,
\end{equation}
\begin{equation}
	E_{i,\alpha}   = \frac{ \sum_{m'=1}^{n_{\mu} }  \mu_{m'}^2\psi_{m',i,\alpha} w_{m'} }{ \sum_{m'=1}^{n_{\mu} }  \psi_{m',i,\alpha} w_{m'} } \, ,
\end{equation}
\begin{equation}
	F_{s,m,i,\alpha}^*  \! =  \! \frac{ \sum_{m'=1}^{n_{\mu} }   \!\!  f_{s,m',m,i} ^* \psi_{m',i,\alpha} w_{m'} }
	{ \sum_{m'=1}^{n_{\mu} }  \!   \psi_{m',i,\alpha} w_{m'} } \,  , \  f_{s,m',m,i} ^*\!  \!=  \! \! f_{s,i}^*(\mu_{m'},\mu_{m} ) \, .
\end{equation}
$\bm{q}_{i}  = [q_{i,L}, q_{i,R}  ]^{\top}$ is the vector of the expansion coefficients of the given source.
The RHS of Eq. \eqref{bte-lld} is defined with
\begin{equation}
	\boldsymbol{\phi}_{i}  = [\phi_{i,L}, \phi_{i,R}]^{\top} \,  , \quad 
	\boldsymbol{J}_{i}  = [J_{i,L}, J_{i,R}]^{\top} \, ,
\end{equation}
which are the vectors of the expansion coefficients of the LOSM solution given by
\begin{equation}\label{dfem-exp-lo}
	\phi_i(x) = \sum_{\alpha=L,R} \phi_{\alpha,i} b_{\alpha,i}(x) \, , \quad
	J_i (x) = \sum_{\alpha=L,R}  J_{\alpha,i} b_{\alpha,i}(x) \, .
\end{equation}
Using these expansions, the LOSM equations are discretized with the LLD Galerkin method.
The discrete LOSM equations are algebraically consistent with the high-order transport equation  
approximated by the LLD scheme given by Eq. \eqref{bte-lld}.

\section{Numerical Results}

We consider transport problems for a slab $x\in [0,1]$ with the following set of different total cross sections and scattering ratios $c=\frac{\sigma_s}{\sigma_t}$:
$\sigma_t=\{ 0.1;1;5;10\}$ and  $c=\{0.5;0.9;1\}$.
The phase function is given by $f_s(\mu_0) = (1+\mu_0)^{10}$ which has the following
coefficients of expansion in Legendre polynomials: 
$f_{s,0}=1$, $f_{s,1}=0.83333$, $f_{s,2}=0.57692$.
There is a constant source of $q=1$. The boundary conditions are vacuum ($\psi_{in}^{\pm}=0$).
The spatial mesh consists of $n_x=10$ intervals. The double $S_8$ Gauss-Legendre quadrature set ($n_{\mu}=16$) is used.  
The convergence criterion with $\varepsilon=10^{-6}$ (see Algorithm \ref{NSM-iter}) is applied.

\begin{figure}[h!]
	\centering
	\begin{minipage}{.5\textwidth}
		\centering
	\includegraphics[scale=0.4]{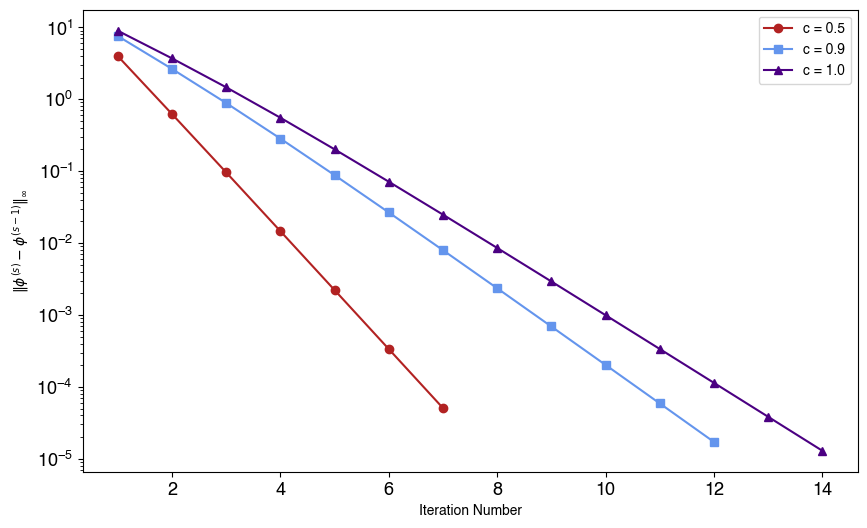}
	\caption{$\|  \!   \phi^{(s)}\!    -  \!\phi^{(s-1)}   \!  \|_{\infty }$ versus iteration number ($s$) for $\sigma_t \! =\! 0.1$}
\label{histories-sgt=0.1}
	\end{minipage}%
	\begin{minipage}{.5\textwidth}
	\centering
\includegraphics[scale=0.4]{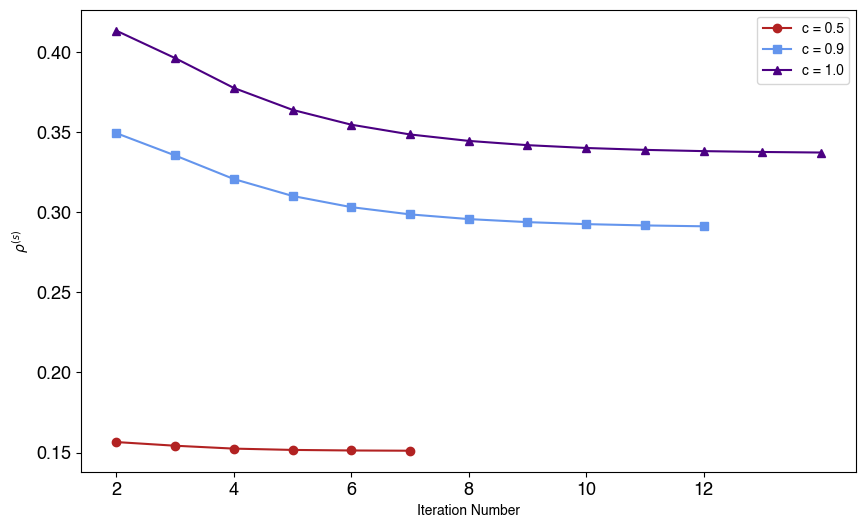}
\caption{$\rho^{(s)}$ versus iteration number ($s$)  for $\sigma_t=0.1$}
\label{rho-sgt=0.1}
	\end{minipage}
\end{figure}

\begin{figure}[h]
	\centering
	\begin{minipage}{.5\textwidth}
		\centering
		\includegraphics[scale=0.4]{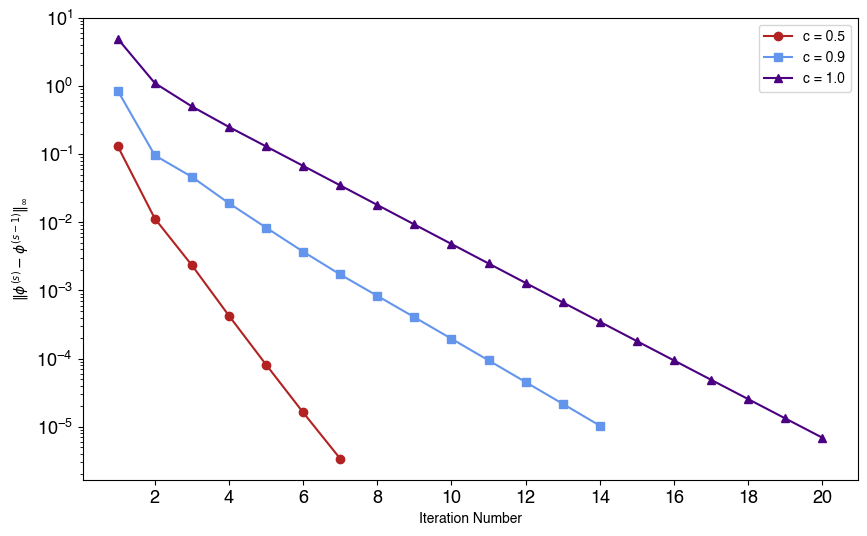}
		\caption{$\|  \!   \phi^{(s)}\!    -  \!\phi^{(s-1)}   \!  \|_{\infty }$ versus iteration number ($s$) for $\sigma_t \! =\! 1$}
		\label{histories-sgt=1}
	\end{minipage}%
	\begin{minipage}{.5\textwidth}
		\centering
		\includegraphics[scale=0.4]{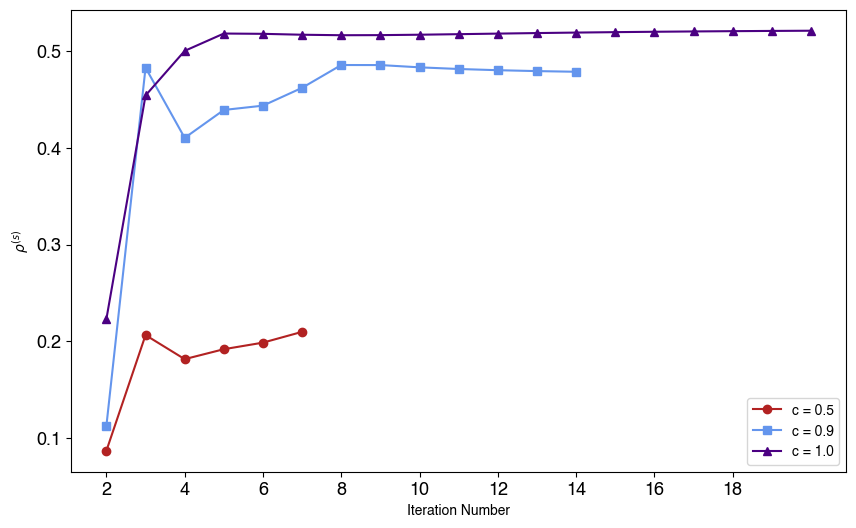}
		\caption{$\rho^{(s)}$ versus iteration number ($s$)  for $\sigma_t=1$}
		\label{rho-sgt=1}
	\end{minipage}
\end{figure}

\begin{figure}[h]
	\centering
	\begin{minipage}{.5\textwidth}
		\centering
		\includegraphics[scale=0.4]{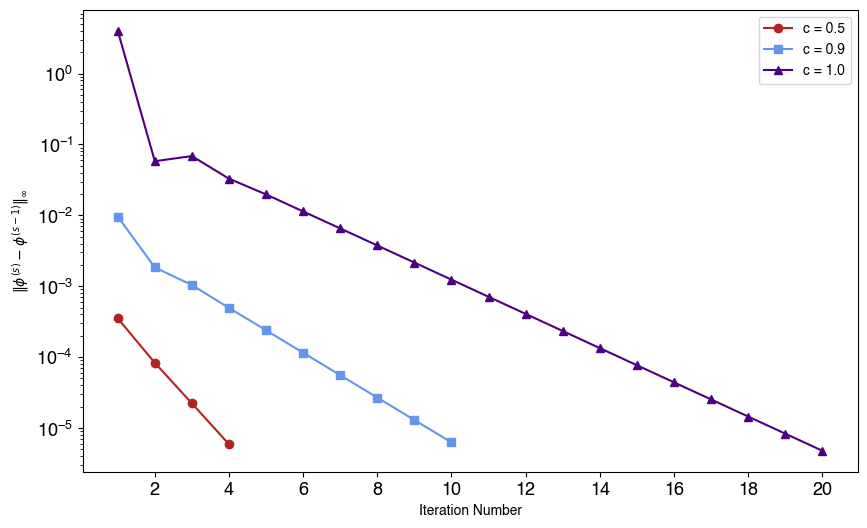}
		\caption{$\|  \!   \phi^{(s)}\!    -  \!\phi^{(s-1)}   \!  \|_{\infty }$ versus iteration number ($s$) for $\sigma_t \! =\! 10$}
		\label{histories-sgt=10}
	\end{minipage}%
	\begin{minipage}{.5\textwidth}
		\centering
		\includegraphics[scale=0.4]{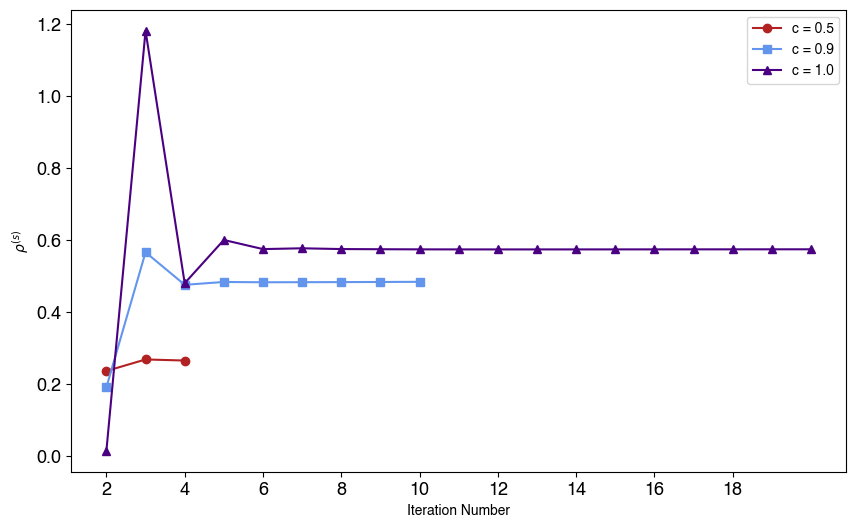}
		\caption{$\rho^{(s)}$ versus iteration number ($s$)  for $\sigma_t=10$}
		\label{rho-sgt=10}
	\end{minipage}
\end{figure}

Table \ref{iterations} presents the number of iterations of  the NSM method in these tests. 
Table \ref{rho} shows the  convergence rate 
$
\rho^{(s)}~=~\frac{ \|  \phi^{(s)} -  \phi^{(s-1)} \|_{\infty } } { \| \phi^{(s-1)} -  \phi^{(s-2)} \|_{\infty } }
$
estimated on the last iteration.
Figures~\ref{histories-sgt=0.1}, \ref{histories-sgt=1}, and 
\ref{histories-sgt=10} present 
$\|  \phi^{(s)} -  \phi^{(s-1)} \|_{\infty}$ versus the iteration number  demonstrating the convergence behavior of the NSM method for  cases with $\sigma_t=0.1;1,10$. 
Figures \ref{rho-sgt=0.1}, \ref{rho-sgt=1}, and 
\ref{rho-sgt=10} show the  convergence rate $\rho^{(s)}$ versus $s$ for these cases.
Table \ref{rho-SI} presents the estimated  rate of convergence of the source iteration method without acceleration \cite{mla-ewl-pne-2002}.

\begin{table}[h!]
	\centering
	\caption{\bf  The NSM Method:  Number of Iterations \label{iterations} }
	\begin{tabular}{|c|c|c|c|c|}
		\hline
		c   \textbackslash   \  \(\sigma_{t}\)\ & 0.1      & 1         & 5        & 10       \\ \hline 
		0.5                            & 7  & 7   & 5 & 4 \\ \hline
		0.9                            & 12 & 14  &  11 & 10 \\ \hline
		1                              & 14    & 20 & 21 & 20 \\ \hline
	\end{tabular}
\end{table}

\begin{table}[h!]
	\centering
	\caption{\bf  The NSM Method:   Estimated Convergence Rates  $\rho$ \label{rho}}
		\begin{tabular}{|c|c|c|c|c|}
			\hline
			c   \textbackslash   \  \(\sigma_{t}\)\ & 0.1      & 1                & 5        & 10       \\ \hline
			0.5                            & 0.151  & 0.209  & 0.243  & 0.265 \\ \hline
			0.9                            & 0.291  & 0.468  & 0.487  & 0.484 \\ \hline
			1                              & 0.371  &  0.521  & 0.571 & 0.574 \\ \hline
		\end{tabular}
\end{table}

    \begin{table}[h!]
	\centering
	\caption{\bf    The SI Method:  Estimated Convergence Rates $\rho$ \label{rho-SI}   }
		\begin{tabular}{|c|c|c|c|c|}
			\hline
			c   \textbackslash   \  \(\sigma_{t}\)\ & 0.1      & 1                & 5        & 10       \\ \hline
			0.5                            & 0.311  & 0.469  & 0.496  & 0.498 \\ \hline
			0.9                            & 0.560  & 0.844  & 0.895  & 0.898 \\ \hline
			1                              & 0.622  &  0.938  & 0.994 & 0.998 \\ \hline
		\end{tabular}
\end{table}

At this stage of research, we have applied the LLD  discretization method to the high-order transport equation. 
It is known to be a more robust scheme compared to the LD Galerkin method.
The LLD and LD schemes are of high-order accuracy, and hence they are non-monotonic and do not preserve positivity of the solution in a general case.  
The prolongation operator involves linear-fractional factors. These factors can be  sensitive to the alternating numerical solution.
This can affect the iteration  behavior.
We plan to apply various monotonization techniques to  preserve the positivity of the numerical solution.


\FloatBarrier

\section{Conclusions}

We have proposed a new nonlinear iteration method for solving particle transport problems with anisotropic scattering based on the SM method.  The high-order transport equation formulated with the nonlinear prolongation operator is discretized  with the LLD Galerkin method. The low-order SM equations are discretized consistently with the high-order transport equation. An independent discretization of the system of  equations of the NSM method can also  be  used because the SM method is known to be stable with this type of discretization  as well.
The linear continuous Galerkin method can be applied to the LOSM equations with the transport equation approximated by the LD scheme
using various mapping techniques \cite{warsa-dya-jctt-2018}.
In future work, we plan  to use the NSM method as a basis for developing a method for transport problems with highly forward-peaked scattering by decomposing the scattering term into smooth and singular components \cite{vyag-aristova-m&c1997,aristova-vyag-2000}.
The formulated method can be naturally extended to multigroup problems.

\section{Acknowledgments}
We wish to thank Anil Prinja for useful discussions.
This work was supported by the Center for \linebreak 
Advancing the Radiation
Resilience of Electronics (CARRE), a PSAAP-IV project funded by
the Department of Energy (DoE) National Nuclear Security Administration (NNSA), award number:  DE-NA0004268.

\bibliographystyle{elsarticle-num}
\bibliography{ija-dya-tans-2026}

\begin{thebibliography}{10}
\expandafter\ifx\csname url\endcsname\relax
  \def\url#1{\texttt{#1}}\fi
\expandafter\ifx\csname urlprefix\endcsname\relax\def\urlprefix{URL }\fi
\expandafter\ifx\csname href\endcsname\relax
  \def\href#1#2{#2} \def\path#1{#1}\fi

\bibitem{dya-vyag-kiam-1985}
D.~Y. Anistratov, V.~Y. Gol'din, Solution of the multigroup transport equation
  by the quasi-diffusion method, Preprint of the Keldysh Institute for Applied
  Mathematics, the USSR Academy of Sciences 128 (1986) 19pp, in Russian.

\bibitem{khattab=larsen-1991}
K.~M. Khattab, E.~W. Larsen, Synthetic acceleration methods for linear
  transport problems with highly anisotropic scattering, Nuclear Science and
  Engineering 107~(3) (1991) 217--227.

\bibitem{morel-manteuffel-1991}
J.~E. Morel, T.~A. Manteuffel, An angular multigrid acceleration technique for
  {$S_n$} equations with highly forward-peaked scattering, Nuclear Science and
  Engineering 107~(4) (1991) 330--342.

\bibitem{vyag-aristova-m&c1997}
V.~Gol'din, E.~Aristova, The method for consideration of a strong scattering
  anisotropy in transport equation, in: Proc. of Joint Int. Conf. on Math.
  Meth. \& Supercomp. in Nucl. App., M\&C 1997, La Grange Park, Il, 1997, pp.
  1507--1516.

\bibitem{pautz-morel-adams-m&c1999}
S.~D. Pautz, J.~E. Morel, M.~L. Adams, Angular multigrid acceleration method
  for {$S_n$} equations with highly forward-peaking scattering, in: Proc. of
  Int. Conf. on Math. and Comp., Reactor Phys. and Envir. Anal. in Nucl. Appl,
  M\&C 1999, Vol.~1, Madrid, Spain, 1999, pp. 647--656.

\bibitem{aristova-vyag-2000}
E.~Aristova, V.~Gol'din, Computation of anisotropy scattering of solar
  radiation in atmosphere (monoenergetic case), Journal of Quantitative
  Spectroscopy and Radiative Transfer 67~(2) (2000) 139--157.

\bibitem{patel-warsa-prinja-2020}
J.~K. Patel, J.~S. Warsa, A.~K. Prinja, Accelerating the solution of the
  {$S_n$} equations with highly anisotropic scattering using the
  {F}okker-{P}lanck approximation, Annals of Nuclear Energy 147 (2020) 107665.

\bibitem{sm-1976}
E.~E. Lewis, J.~W.~F.~Miller, A comparison of {$P_1$ }synthetic acceleration
  techniques, Transactions of the American Nuclear Society 23 (1976) 202--203.

\bibitem{mla-ewl-pne-2002}
M.~L. Adams, E.~W. Larsen, Fast iterative methods for discrete-ordinates
  particle transport calculations, Progress in Nuclear Energy 40 (2002) 3--159.

\bibitem{warsa-dya-jctt-2018}
J.~S. Warsa, D.~Y. Anistratov, Two-level transport methods with independent
  discretization, Journal of Computational and Theoretical Transport 47~(4-6)
  (2018) 424--450.

\end{thebibliography}

\end{document}